\documentstyle{amsppt}
\define\2x2#1#2#3#4{\pmatrix #1&#2\\#3&#4\endpmatrix}
\define\M#1#2{M^{#1#2}}
\define\m#1#2{m^{#1#2}}
\define\eps{\varepsilon}
\define\g{\goth}
\define\wt{\widetilde}
\define\A#1{\Cal A(#1)}
\define\F#1{\Cal F(#1)}
\define\N#1{\Cal N(#1)}
\define\I#1{\Cal I(#1)}
\define\Si#1{\Cal S(#1)}
\define\T#1{#1\widehat\otimes#1^*}
\define\Bn#1{\Cal B(#1)}
\define\Aa#1#2{\Cal A(#1,#2)}
\define\Ff#1#2{\Cal F(#1,#2)}
\define\Ii#1#2{\Cal I(#1,#2)}
\define\Nn#1#2{\Cal N(#1,#2)}
\define\Bb#1#2{\Cal B(#1,#2)}
\define\ot#1{\underset#1\to{\otimes}}
\define\pot#1{\underset#1\to{\widehat\otimes}}
\define\wot#1{\underset{\bold#1}\to{\widetilde\otimes}}
\define\simm{\underset M\to\sim}
\define\tens#1{#1\pot{#1}#1}
\define\Ho#1#2{\Cal H^#1(#2,#2^*)}
\define\apb{{}_{\g A}P_{\g B}}
\define\bqa{{}_{\g B}Q_{\g A}}
\define\ap{{}_{\g A}P}
\define\bq{{}_{\g B}Q}
\define\lp#1{\ell_p(#1)}
\define\LP#1{L_p(\mu,#1)}
\define\nrm#1{\Vert #1\Vert}
\define\jnrm#1{\Vert #1\Vert_{\g J_p}}
\define\jnrmx#1{\Vert #1\Vert_{\g J_p(X)}}
\define\pb{P_{\g B}}
\define\qa{Q_{\g A}}
\define\LIM{\operatorname{LIM}}
\define\tr{\operatorname{tr}}
\define\Tr{\operatorname{Tr}}
\define\rk{\operatorname{rank}}
\define\cl{\operatorname{cl}}

\define\clspan{\operatorname{clspan}}
\define\way{weak a\-men\-a\-bi\-li\-ty}
\define\wae{weak\-ly a\-men\-ab\-le}
\redefine\l{\langle}
\define\r{\rangle}
\define\set#1{\{#1\}}
\refstyle{A}
\topmatter
\title Factorization and weak amenability of $\A X$\endtitle
\author Niels Gr\o nb\ae k \endauthor
\address Department of Mathematics, Institute for Mathematical
Sciences, Universitets\-parken 5, DK-2100
Copenhagen {\O}, Denmark\endaddress
\email gronbaek\@math.ku.dk\endemail
\abstract We investigate weak amenability of the Banach algebra $\A X$ of
approximable operators on a Banach space $X$ and its relation to
factorization properties of operators in $\A X$. We show that if
$\A X$ is weakly amenable, then either $\A X$ is self-induced (a nice
factorization property), or $X$ is very special, combining some of the
exotic properties of
the spaces of Gowers and Maurey [GM] and of Pisier [P1]. In the class of self-induced Banach
algebras we show that weak amenability is preserved under
an equivalence of Morita type. Using this we extend some results of
A\. Blanco [B1, B2] about weak amenability of $\A X$. \endabstract 
\endtopmatter
\document

\subheading{0. Introduction}

Recall that a Banach algebra $\frak A$ is called {\it \wae} if every
bounded derivation $D\colon\frak A\to\frak A^*$ is inner, or
equivalently if the first bounded Hochschild cohomology group
$\Ho1{\frak A}$ vanishes. Recently \way\ has been investigated for
algebras of the type $\A X$ for an infinite dimensional Banach space
$X$. In [DGG] it is shown that $\A X$ is \way\, when $X=\ell_p(Y)$ with
$Y$ reflexive and having the approximation property, or when
$X=E\oplus C_p$, where $E$ has the bounded approximation property and
$C_p$ denotes any of the universal spaces introduced by W.B\. Johnson
in [Jo]. In \cite{B2} Blanco introduces a technical property of $X$,
socalled trace unbounded triples that allows for taking averages of
matrix-like representations of a given finite rank operator. Using
this Blanco establishes \way of $\A X$ for a wide range of Banach
spaces $X$. In \cite{B1} Blanco studies hereditary properties of as
well as necessary conditions for \way\ for algebras $\A X$. In this,
factorization properties play a crucial role.

In the present paper we shall take a approach almost exclusively
related to factorization properties. We show that if $\A X$ is \wae,
then either $X$ is pathological (probably non-existing), or $\A X$ is
socalled self-induced. Self-induced Banach algebras constitute the
class of Banach algebas for which a Morita theory can naturally be
develloped. Hence our approach will be to tranfer Hochschild
cohomology from a few key examples by means of Morita equivalence,
i.e\. by means of factorization properties. In this way we give a
unified approach to some of Blanco's results with shorter and less
technical proofs and in some cases improvements of the statements. We
emphasise though, that our aim is to view the question of \way\ from a
more general stand i.e\.  that of Morita equivalence. Blanco's results
serve here as a source of test cases. In particular, the result
\cite{B2} on \way\ of $\A T$ with $T$ the Tsirelson space, remains a
challenge.

\subheading{1. Preliminaries} For  Banach spaces $X$ and $Y$ we
consider the following spaces of operators 
$$
\align
&\Ff XY=\{\text{finite rank operators  $X\to Y$}\}\\
&\Aa XY=\{\text{approximable operators  $X\to Y$}\}\\ 
&\Nn XY=\{\text{nuclear operators  $X\to Y$}\}\\ 
&\Ii XY=\{\text{integral operators  $X\to Y$}\}\\
&\Bb XY=\{\text{bounded operators  $X\to Y$}\}.
\endalign
$$
As it is customary, we shall write $\operatorname{Operators}(X)$ for
$\operatorname{Operators}(X,X)$. We write $|\cdot|_{\Cal
N},|\cdot|_{\Cal I}$, and $\nrm\cdot$ for the nuclear, integral, and
uniform norm on $\Nn XY,\Ii XY$, and $\Bb XY$ respectively. The identity operator on $X$ is
denoted by $\bold 1_X$, or if the context is clear, simply by $\bold1$. 

\smallskip
For Banach spaces $E$ and $F$ their projective tensor product is
denoted $E\pot{}F$. The {\it tensor algebra} of $X$ is $\T X$ with
multiplication given by 
$$
(x\otimes x^*)(\xi\otimes \xi^*)=x^*(\xi) x\otimes\xi^*,\quad
x,\xi\in X;\;x^*,\xi^*\in X^*.
$$
The trace $\tr\colon
X\pot{}X^*\to\Bbb C$ and operator trace $\Tr\colon
X\pot{}X^*\to\N X$ are given by 
$$
\tr(x\otimes x^*)=x^*(x),\;\Tr(x\otimes x^*)(\xi)=x^*(\xi)x,\quad
x,\xi\in X,x^*,\xi^*\in X^*
$$
Note that $\Tr$ maps onto $\N X$. 

For Banach spaces $X$ and $Y$ we denote the statement '$X$ is
isomorphic to $Y$' by $X\cong Y$. For spaces in duality we shall use
$\l\cdot\,,\cdot\r$ to denote the corresponding bilinear form, in
particular we shall write $\l x,x^*\r=x^*(x)$ for $x\in X,x^*\in
X^*$. We note in particular the {\it trace duality}
$$
\l F,T\r=\tr(F^*T),\quad F\in\F X,T\in\Bn{X^*},
$$
which isometrically identifies $(\T X)^*$ and $\A X^*$ with $(\Bn{X^*},\nrm\cdot)$
and $(\I{X^*},|\cdot|_{\Cal I}$ respectively.

For any normed space $E$ the unit ball is denoted by $E_1$.

Let $X_n,\;n\in\Bbb N$ be a sequence of Banach spaces. We denote the
$\ell_p$-sums of this sequence by $(\oplus_1^\infty X_n)_p$ for
$p=0,1\leq p\leq\infty$. If $X_n=X$ for all $n\in\Bbb N$ we simply
write $c_0(X)$, or $\ell_p(X),\;1\leq p\leq\infty$.

We shall frequently without further reference use the fact (see [D]) that for a
Banach space $X$
$$
\multline
\A X\text{ has a bounded left approximate identity}\iff\\ X\text{ has
the bounded approximation property}.\endmultline
$$ 
Hence, if $X$ has the bounded approximation property, we may use Cohen
factorization in the Banach algebra $\A X$.
\smallskip
The definitions of  Banach (co)homological concepts are standard and
can be found for example in [H] and [J]. We shall only here point to

\definition{1.1 Definition} Let $\g A$ be a Banach algebra, let $X$ be a right
Banach $\g A$-module, and let $Y$ be a left Banach $\g A$-module. We
define
$$
X\pot{\g A}Y=X\pot{}Y/N,
$$
where $\pot{}$ is the projective tensor product and
$N=\clspan\{x.a\otimes y-x\otimes a.y\mid x\in X,y\in Y,a\in \g 
A\}$. Thus, $X\pot{\g A}Y$ is the universal object for linearizing
bounded, $\g A$-balanced, bilinear maps $X\times Y\to Z$.
\enddefinition
\medskip
We start by recalling some facts about bounded derivations $D\colon\A
X\to\A X^*$. As said above, we identify $\A X^*$ with $\I{X^*}$,
and $(X\pot{}X^*)^*$ with $\Bn{X^*}$ via trace duality. Consider the
diagram 
$$
\CD
X\pot{}X^*      @>>\Tr>        \A X\\
@V\delta VV                      @VVDV\\
\Bn{X^*}        @<\Tr^*<<      \I{X^*}
\endCD
$$
Since $X\pot{}X^*$ is biprojective [S] and in particular \wae, the
derivation $\delta$ is inner. Hence we have
\proclaim{1.2 Proposition} Let $D\colon\A
X\to\A X^*$ be a bounded derivation. Corresponding to $D$
there is $T\in\Bn{X^*}$ such that 
$$
\l F,D(G)\r=\tr((FG-GF)^*T),\quad F,G\in\F X.
$$ 
Consequently $D$ is inner if and only if  $T\in\I{X^*}+\Bbb
C\bold1_{X^*}$, that is, if and only if there is $\lambda\in\Bbb
C$ and $K>0$ such that 
$$
\vert\tr( F^*T-\lambda F^*)\vert\leq K\nrm F,\quad F\in\F X.
$$
\endproclaim

\subheading{2. Factorization properties and \way}

An important aspect of Morita theory is to provide tools to compare
homological properties of Banach algebras using \lq good factorization
properties'. In this section we shall extract such factorization
properties in order to compare $\Ho n{\g A}$ and $\Ho n{\g B}$ for
Banach algebras $\g A$ and $\g B$. Our focus shall be on $n=1$ and
Banach algebras of the type $\A X$.

First we make precise what is meant by \lq good factorization':

\definition{2.1 Definition} A Banach algebra $\g A$ is called {\it
  self-induced} if 
$$
\g A\cong\g A\pot{\g A}\g A
$$
The Banach algebra $\A X$ {\it factors approximately through} $Y$, if
$$
\A X\cong\Aa YX\pot{\A Y}\Aa XY,
$$
where the isomorphisms are implemented by multiplication.
\enddefinition

The usefulness of these factorization properties is that one may
define linear maps in terms of balanced bilinear maps. A key example is

\proclaim{2.2 Lemma} Suppose that $\g A$ is self-induced. Let $D\colon
\g A\to \g A^*$ be a bounded derivation. Let $\g B$ be a Banach
algebra which contains $\g A$ as a closed 2-sided ideal. Then $D$
may be extended to a bounded derivation $\widetilde D\colon \g B\to\g
A^*$. 
\endproclaim
\demo{Proof} Let $T\in \g B$ and consider the bilinear map 
$\Phi_T\colon \g A\times\g A\to\Bbb C$ given by
$$
\Phi_T(a,b)=\l a,D(bT)\r-\l Ta,D(b)\r,\quad a,b\in\g A.
$$
Then $\Phi_T$ is balanced, i.e. $\Phi_T(ac,b)=\Phi_T(a,cb)$, so we may
define
$$
\l ab,\widetilde D(T)\r=\Phi_T(a,b),\quad a,b\in\g A.
$$
One checks  that this defines a bounded derivation $\g B\to\g
A^*$ extending $D$ (in fact the only possible such).
\enddemo
\remark{Remark} In the same way $D$ can be lifted to a derivation $\g
A\to\g B^*$.
\endremark

\example{2.3 Example} Assume that the
multiplication
$$
\A X\pot{}\A X\to\A X
$$
is surjective. If $X$ in addition has the approximation property, then
$\A X$ is self-induced. Suppose namely 
$$
\sum A_nB_n=0
$$ 
with $A_n@>>> 0$ and $\sum \Vert
B_n\Vert<\infty$, and let $\eps>0$. Since $X$ has the approximation
property, we may choose 
$U\in\A X$ so that 
$$
\sup_n\Vert UA_n-A_n\Vert\leq\eps.
$$
Then 
$$
\align
\Vert\sum A_n\ot{\A X}B_n\Vert&\leq\Vert\sum UA_n\ot{\A
  X}B_n\Vert+\sum\Vert A_n-UA_n\Vert\Vert B_n\Vert\\
&\leq 0+\eps\sum\Vert B_n\Vert.
\endalign
$$

Thus, an important case occurs, when $X$ has the bounded approximation
property, using the bounded approximate identity in $\A X$ ([D]). 
\endexample

\smallskip 
The approximation property is not essential here as will be clear in
the course of the paper. However, in the case of nuclear operators
self-inducedness and the approximation property is one and the same thing.

\proclaim{2.5 Proposition} The Banach algebra $\N X$ is self-induced
if and only if $X$ has the approximation property.
\endproclaim
\demo{Proof} We define a bounded balanced bilinear form on $\N X$ by
$$
\phi(N,M)=\tr(UV)
$$
where $U,V\in X\pot{}X^*$ with $\Tr(U)=N$ and $\Tr(V)=M$. This is well defined, since if
$\Tr(U)=0$, then $U(X\pot{}X^*)=(X\pot{}X^*)U=\{0\}$. Suppose now that
$\N X$ is self-induced. Then $\phi$ defines a bounded linear
functional on $\N X$ which agrees with the standard trace on $\F
X$. But then $X$ must have the approximation property. Conversely, if
$X$ has the approximation property, then $\N X$ is isometrically
isomorphic to the tensor algebra $\T X$. Since any rank one tensor
$x\otimes x^*$ has a factorization $p(x\otimes x^*)$ with $p$ a norm 1
projection, it follows easily that $\T X$ is self-induced.
\enddemo

In order to investigate the relation between derivations and
self-inducedness we look at derivations of the following type.

Let $\phi\in(\g A\pot{\g A}\g A)^*$. For convenience we  shall use the same
symbol $\phi$ for the corresponding balanced bilinear functional.
Then one checks that  
$$
  \l a,D(b)\r:=\phi(a,b)-\phi(b,a)
$$
defines a bounded derivation $D\colon\g A\to\g A^*$. 

How does this look in the setting of $\g A=\A X$? First we need to describe
balanced bilinear functionals.

\proclaim{2.6 Lemma} Let $\phi\in(\Aa YX\pot{\A Y}\Aa XY)^*$. Then there
is $T\in\Bn{X^*}$ such that 
$$
\phi(F,G)=\tr((FG)^*T)\quad F\in\Ff YX,G\in\Ff XY
$$
and
$$
\sup\{\vert\tr((FG)^*T)\vert\mid F\in\Ff YX _1,G\in\Ff XY_1\}<\infty.
$$

If there is a constant $K>0$ such that every $F\in \F
X$ is a finite sum of products, $F=\sum_0^n U_iV_i$ with $U_i\in\Ff YX,V_i\in\Ff
XY$ and $\sum_0^n\Vert 
U_i\Vert\Vert V_i\Vert\leq K\Vert F\Vert$, then $\Aa YX\pot{\A Y}\Aa 
XY\cong \A X$.
\endproclaim

\demo{Proof} To find $T\in \Bn{X^*}$ we apply $\phi$ to rank-1
operators. Choose $e\in X,\;e*\in X^*$ with $\l e,e^*\r=1$. Define
$T\in\Bn{X^*}$ by
$$
\l x,T(x^*)\r=\phi(x\otimes e^*,e\otimes x^*),\;x\in X,x^*\in X^*.
$$ 
Then
$$
\align
\phi(x\otimes x^*,y\otimes y^*)&=\phi((x\otimes e^*)(e\otimes
x^*),y\otimes y^*)\\
&=\phi(x\otimes e^*,(e\otimes x^*)(y\otimes y^*))\\
&=\l y,x^*\r\phi(x\otimes e^*,e\otimes y^*)\\
&=\l y,x^*\r\l x,T(y^*)\r\\
&=\tr([(x\otimes x^*)(y\otimes y^*)]^*T).
\endalign
$$
By linearity $\phi$ is given by $T$ as claimed. The norm estimate is just
$$
\Vert\phi\Vert=\sup\{\vert\tr((FG)^*T)\vert\mid F\in\Ff YX _1,G\in\Ff
XY_1\}<\infty. 
$$

Now assume that we have factorization with $K>0$ as described. Clearly
then the multiplication 
$$
\mu\colon\Aa YX\pot{\A Y}\Aa XY\to \A X
$$
is surjective. Let $\phi\in(\Aa YX\pot{\A Y}\Aa XY)^*$ with corresponding
$T\in\Bn{X^*}$.
The assumption  ensures that
$$
f(\sum U_iV_i)=\tr((\sum U_iV_i)^*T)
$$
defines a bounded functional, so that the dual map 
$$\mu^*\colon\A
X^* \to (\Aa YX\pot{\A Y}\Aa XY)^*
$$ 
is also surjective, i.e. $\mu$ is
an isomorphism.
\enddemo

In Section 4 of [B1] Blanco discusses necessary conditions for \way. He
shows that if $\A X$ is \wae, then either $X$ is indecomposable (i.e\.
$X$ is not the direct sum of two infinite-dimensional Banach spaces) or the
trace defines an unbounded bilinear map associated with a
decomposition. These considerations are naturally futher explored by  means
of self-inducedness. The Banach spaces of Pisier [P1] for which $\A X=\N
X$ are crucial in this. We note some simple reformulations of
this property. But first we need the following estimate of norms,
which is essentially an elaboration of the proof of [DU, Theorem
VIII.4.12]
\proclaim{2.7 Lemma} Let $T\in\Ii XY,S\in \Bb YZ$. If $S$ is weakly
compact, then 
$$
|ST|_{\Cal N}\leq\nrm S|T|_{\Cal I}
$$
\endproclaim
\demo{Proof} That $ST$ is nuclear is the statement of [DU, Theorem
VIII.4.12.(i)]. We note from the proof of this, that 
$S$ being weakly compact, there is a reflexive space $W$ and operators
$A\in\Bb YW$ and $B\in\Bb WZ$ such that $S=BA$ by [DU, Corollary
VIII.4.9]. Furthermore, a close inspection of the proof shows that
$\nrm S=\inf\set{\nrm B\nrm A}$, where the infimum is taken over such
factorizations. Reasoning along with [DU] we get
$$
|AT|_{\Cal N}=|AT|_{\Cal I}\leq\nrm A|T|_{\Cal I},
$$ 
so that 
$$
|ST|_{\Cal N}=|BAT|_{\Cal N}\leq \nrm B|AT|_{\Cal N}\leq\nrm B\nrm
A|T|_{\Cal I}.
$$
Taking the infimum over $\nrm A\nrm B$ gives the wanted estimate.
\enddemo  

\proclaim{2.8 Proposition} Let $X$ be a Banach space. Then
(i)$\implies$(ii)$\iff$(iii)$\implies$(iv), where
\roster
\item"(i)" $\A X=\N X$
\item"(ii)" $\A X\subseteq \I X$
\item"(iii)" There is $C>0$ such that 
$$
|\tr(AB)|\leq C\Vert A\Vert\Vert B\Vert,\quad A,B\in \F X
$$
\item"(iv)" the multiplication $\A X\pot{}\A X\to\A X$ maps onto $\N X$.
\endroster

In particular, if the multiplication $\A X\pot{}\A X\to\A
X$ is known to be surjective, then all four are equivalent.
\endproclaim
\demo{Proof} Since in general $\N X\subseteq\I X$, (i)$\implies$(ii)
is obvious. 
\flushpar (ii)$\implies$(iii): An  application of the
closed graph theorem shows that the inclusion $(\A
X,\Vert\,\cdot\,\Vert)\hookrightarrow(\I X,\vert\,\cdot\,\vert_\Cal I)$ is
continuous, thus providing $C>0$ so that $|A|_\Cal I\leq C\Vert
A\Vert$ and hence $|A^*|_\Cal I\leq C\Vert A\Vert$
for all $A\in\A X$, since $|A|_\Cal I=|A^*|_\Cal I$ ([DJT, Theorem 5.15]). This gives
$$
|\tr(AB)|=|\tr(B^*A^*)\leq|A^*|_\Cal I\Vert B\Vert\leq C\Vert
A\Vert\Vert B\Vert 
$$
for all $A,B\in\F X$.
\flushpar(iii)$\implies$(ii): By one of the definitions [DJT,??] of the
integral norm (iii) states that $|A^*|_\Cal I\leq C\Vert A\Vert$ for
each $A\in\F X$ from which (ii) follows, again using $|A|_\Cal I=|A^*|_\Cal
I$.
\flushpar(iii)$\implies$(iv): By Lemma 2.7  we have for $A,B\in\F X$
$$
|AB|_\Cal N\leq\Vert A\Vert|B|_\Cal I\leq C\Vert A\Vert\Vert B\Vert
$$
Hence if $A=\sum A_nB_n$ with $A_n,B_n\in\A X,\;\sum \Vert
A_n\Vert\Vert B_n\Vert<\infty$ we have that 
$$
\sum|A_nB_n|_\Cal N\leq C\sum \Vert A_n\Vert\Vert B_n\Vert<\infty
$$
so that the series is absolutely convergent in the nuclear norm and
thus $A\in\N X$. Since the
multiplication $\N X\pot{}\N X\to\N X$ always is surjective, we arrive
at (iv).
\enddemo

We shall now show that weak amenability of $\A X$ forces either $\A
X$ to be self-induced or the underlying space $X$ to be very peculiar,
combining some of the pathological properties of the spaces of Pisier
[P1] and Gowers and Maurey [GM].

\proclaim{2.9 Theorem} Suppose that  $\A X$ is not self-induced. Then
$\A X$ is weakly amen\-able if and only if both (a) and (b) hold, where 
\roster
\item"(a)" $\A X=\N X$
\item"(b)" The kernel, $K,$ of the operator trace $\Tr\colon X\pot{} X^*\to \N
X$ is 1-dimensional.
\endroster
If (a) and (b) hold then
\roster
\item"(c)" We have 
\comment
\endcomment
$$
\Bn X=\I X\oplus\Bbb C\bold 1_X\text{ and }\Bn{X^*}=\I{X^*}\oplus\Bbb
C\bold 1_{X^*}.
$$
\item"(d)" If $X=Y\oplus Z$, then either $Y$ or $Z$ is finite
dimensional (i.e\. $X$ is indecomposable). Similarly, $X^*$ is indecomposable.
\endroster
\endproclaim
\demo{Proof} Suppose that $\A X$ is weakly amenable. First we note
that multiplication is surjective, since the map $\A X\pot{}\A X\to\A
X\colon F\otimes G\mapsto FG-GF$ has closed range, and $\A X$ has no
bounded traces. Hence all four conditions of Proposition 2.8 are
equivalent. Now let $\phi\in(\A
X\pot{\A X}\A X)^*$ and let $T_\phi\in\Bn {X^*}$ be the corresponding
linear operator according to Lemma 2.6. Since $\A X$ is weakly amenable it
follows from the paragraph preceeding the lemma, that there is an
integral operator $T\in\I{X^*}$ and $\lambda\in\Bbb C$ so that 
$$
T_\phi=T+\lambda\bold 1_{X^*}.
$$
If $\A X$ is not self-induced we may choose $\phi$ so that
$\lambda\neq 0$. In this case we have
$$
|\tr(AB)|=|\tr(AB)^*|\leq C\Vert A\Vert\Vert B\Vert, \quad A,B\in\F X
$$
for appropriate $C>0$. The statement (a) now follows from
Proposition 2.7(i). 

In general $\N X$ is weakly amenable if and only if
$\dim K\leq 1$ ([G1]). Noting that $X$ does not have the approximation
property (if it were so, $\A X$ 
would be self-induced, cf\. Example 2.3) we arrive at (b). (Recall
that $K=\{0\}\iff X$ has the approximation property.)  

Setting $K=\Bbb Cu$ with $u=\sum x_n\otimes x_n^*$ and $\tr
u=\sum <x_n,x_n^*>=1$ the functional $\varphi\colon\Bn {X^*}\to\Bbb C$
given by
$$
\varphi(T)=\sum <x_n,T(x_n^*)>, \quad T\in\Bn {X^*}
$$
is multiplicative with kernel $\I {X^*}$ (see the proof of Corollary 4
of [G1]). Using the trace duality between $X\pot{}X^*$ and $\Bn{X^*}$
we find that
$$
(\N X)^*=K^\perp=\ker\varphi\tag*
$$
Hence, when (a) holds, we get $\ker\varphi=\I{X^*}$, thus proving the
last equality in (c). The first equality follows by means of the multiplicative
linear functional $T\mapsto\varphi(T^*)$ on $\Bn X$ and the fact that
$T$ is integral if and only if $T^*$ is integral [DJT, Theorem 5.15].

To prove (d) first note that if $P\in\Bn{X^*}$ is a projection then
$\varphi(P)$ is either 1 or 0. A simple applications of the closed
graph theorem gives that the integral and uniform norms are
equivalent on $\I{X^*}$. Accordingly there is a constant $C>0 $ so that 
$$
|\tr(A^*T)-\varphi(T)\tr A|\leq C\Vert A\Vert,\quad A\in\F X,\;T\in
\Bn{X^*}.\tag** 
$$
If $P\in\Bn X$ is a projection with $\rk P=\infty$, we may for each
$n\in\Bbb N$ choose a projection $Q_n\in\F X$ with
$$
\tr Q_n=n\,,\Vert Q_n\Vert\leq n^{\frac12},\text{ and } PQ_n=Q_n\;,
$$
cf\. [P1, Theorem 1.14]. If this goes along with (**) we must have $\varphi(P^*)=1$. Thus, if
$X$ were decomposable, we would have $2=1$. If $P'\in\Bn{X^*}$
is a projection of infinite rank, choose projections $Q_n'\in\F{X^*}$ as above. We may use local reflexivity to
modify the $Q_n'$ to obtain  projections $Q_{\!*\,n}\in\F X$ such that
$$
\tr Q_{\!*\,n}=n\,,\Vert Q_{\!*\,n}\Vert\leq 2n^{\frac12},\text{ and } P'Q_{\!*\,n}^*=Q_{\!*\,n}^*\;,
$$
so that also $X^*$ is indecomposable. 
\enddemo

\subheading{3. Weak amenability of self-induced Banach algebras} From now
on we shall concentrate on 
self-induced Banach algebras.  In order to compare cohomology of such we
shall exploit the double complex 
of Waldhausen [DI]. First we consider the lower
left hand corner of a general  double co-complex in the first quadrant:
$$
\CD
0@>>>\M03@>>>\\
@.    @AAA     @AAA\\
0@>>>\M02@>>>\M12@>>>\\
@.    @AAA        @AAA   @AAA \\
0@>>>\M01@>>>\M11@>>>\M21@>>>\\
@.  @AAA        @AAA @AAA        @AAA\\
@.0@>>>\M10@>>>\M20@>>>\M30\\
@.@.  @AAA        @AAA @AAA\\
@.@.0@.0@.0
\endCD\tag$\Cal D$
$$
The upper indices are meant as coordinates in the first quadrant. We
shall assume that the diagram is commutative.  On the horizontal axis we define the
cohomology $\Cal H^n_h$ 
as kernel modulo image of
$$
@>>>\M0{n+1}@>>>.
$$

The cohomology on the vertical axis, $\Cal H^n_v$, is defined
analogously. We want to compare $\Cal H^1_h$ and $\Cal H^1_v$. In
essense this consists of showing that the associated spectral sequence
collapses at appropriate $E^2$-terms. However, we give a direct
construction of a comparing map using an ad hoc diagram chase.
 
\proclaim{3.1 Lemma} Consider the diagram ($\Cal D$). If there is vertical
exactness at coordinates (1,1), (2,0), (3,0), (1,2), and (2,1), then
we may define a linear map 
$$
\goth D\colon\Cal H^1_v\to\Cal H^1_h
$$
such that 
\roster
 \item"" If there is horizontal exactness at (1,1) and (0,2), then
$\goth D$ is injective. 
\item"" If there is horizontal exactness at (1,2), (2,1), and (0,3),
then $\goth D$ is surjective.
\endroster
\endproclaim

\demo{Proof} First we describe a procedure to associate a cocycle at
$(2,0)$ to each cocycle at $(0,2)$. We adopt the convention
that indices on cochains indicate belonging, i.e\. $\m ij\in\M
ij\,,\mu^{ij}\in\M ij$ ect\. Let $\m02$ be a vertical
cocycle. The numbers at the arrows show the progression in the
diagram chase.
$$
\CD
0@>2.>>0\\
@A1.AA     @A4.AA\\
\m02@>3.>>\m12@>6.>>0\\
@.        @A5.AA   @A8.AA \\
@.\m11@>7.>>\m21@>10.>>0\\
@.  @.   @A9.AA        @A12.AA\\
@.@.\m20@>11.>>\m30\\
\endCD
$$
The existence of $\m11$ and $\m20$ is due to vertical exactness. Since
the vertical map at (3,0) is injective, $\m30=0$ so that $\m20$ is a
cocycle. 

 We next show that this actually gives a map into $\Cal H^1_h$. Hence
suppose that for the same $\m02$ we have made  different choices
$\m11_*,\m20_*$ . Then
$$
\CD
0\\
@A1.AA\\
\m 11-\m11_*@>3.>>\m21-\m21_*\\
@A2.AA@A5.AA\\
\mu^{10}@>4.>>\mu^{20}
\endCD
$$
Here $\mu^{10}$ exists by vertical exactness at (1,1). Since we have
vertical injectivity at (2,0), we must have $\mu^{20}=\m20-\m20_*$, so
$\m20$ and $\m20_*$ are cohomologous.

We next show that this map lifts to the desired map $\goth D$. Hence
assume that $\m02$ cobounds vertically. Then the procedure gives
$$
\CD
\m02@>2.>>\m12\\
@A1.AA@A4.AA\\
\m01@>3.>>\m11@>5.>>0\\
@.@.@A6.AA\\
@.@.0
\endCD,
$$
i.e. coboundaries go to coboundaries.

Now assume that there is horizontal exactness at places (1,1) and
(0,2) and that the procedure $\m02\mapsto\m20$ has resulted in
a coboundary. This is described in
$$
\CD
\m11_*@>4.>>\m21\\
@A3.AA@A1.AA\\
\m10_*@>2.>>\m20\\
\endCD
$$
Elements from the procedure are unstarred $\m ij$'s. We next get
$$
\CD
\mu^{02}@>5.>>\m12\\
@A4.AA@A1.AA\\
\mu^{01}@>3.>>\m11-\m11_*@>2.>>0
\endCD
$$
The arrow 1\. is valid because $\m11_*$ is a vertical coboundary, and
$\mu^{01}$ exists by horizontal exactness at (1,1). By horizontal
injectivity at (0,2) we must have $\mu^{02}=\m02$, i.e\. $\m02$ is a
coboundary. 

Finally assume horizontal exactness at (2,1) and (1,2) and let
$\m20$ be a horizontal cocycle. Then we get the diagram
$$
\CD
\m03 @>12.>> 0\\
@A11.AA @A10.AA\\
\m02 @>9.>> \m12 @>8.>> 0\\
@.@A7.AA @A6.AA\\
@.\m11@>5.>>\m21@>4.>>0\\
@.@. @A3.AA @A2.AA\\
@.@.\m20 @>1.>>0
\endCD
$$
Here $\m11$ and $\m02$ exist due to horizontal exactness at (2,1) and
(1,2). Since we have horizontal injectivity at (0,3) we get $\m03=0$,
altogether showing that the found $\m02$ is a cocycle and by the
procedure is taken to the given $\m20$, i.e\. $\frak D$ is surjective.
\enddemo

\remark{3.2 Remark} Note that in order to define the map $\g D$ we did
not use the assumptions of vertical exactness at plase (1,2) in
full. All that is needed, is that the cocycle $\m12$ corresponding to
the cocycle $\m02$ is actually a coboundary. 
\endremark  
\medskip
We want to use Lemma 3.1 to establish instances of Morita invariance of
Hoch\-schild cohomology, essentially by refining the arguments in
[G2]. The definition of Morita equivalence is usually given in terms of
functors between categories of modules. We only need a slightly weaker
concept (which in the case of Banach algebras with bounded one-sided
approximate identities coincides with the full version of Morita equivalence). 

\definition{3.3 Definition} Let $\g A$ and $\g B$ be self-induced Banach
algebras. Then $\g A$ and $\g B$ are {\it M-equivalent}, in
symbols $\g A\simm\g B$, if there are bimodules $\apb$ and $\bqa$ and 
balanced pairings
$$
[\cdot]\colon\ap\pot{\g B}\qa\to\g A\text{ and
}[\cdot]\colon\bq\pot{\g A}\pb\to\g B\;, 
$$
which are bimodule isomorphisms satisfying
$$
[p\pot{\g B}q].p'=p.[q\pot{\g A}p']\text{ and }[q\pot{\g
A}p].q'=q.[p\pot{\g B}q']\;,\quad p,p'\in \apb,\;q,q'\in \bqa\,.
$$
\enddefinition

 The double complex to which we shall apply the
Theorem 3.1, is the dual complex of the Waldhausen double complex [DI]. We shall for
short write $P$ and $Q$ instead of $\apb$ and $\bqa$. The lower left
hand corner of the Waldhausen double complex is
$$
\CD
0 @<<< \g B\pot{}\g B\pot{}\g B @<<<\\
@.    @VVV     @VVV\\
0@<<<\g B\pot{}\g B@<<< P\pot{}\g B\pot{}Q @<<<\\
@.    @VVV        @VVV   @VVV \\
0@<<< \g B @<<< P\pot{} Q @<<<P\pot{} Q \pot{}\g A @<<<\\
@.  @VVV        @VVV @VVV        @VVV\\
@.0@<<<\g A @<<<\g A\pot{}\g A@<<<\g A\pot{}\g A\pot{}\g A\\
@.@.  @VVV        @VVV @VVV\\
@.@.0@.0@.0
\endCD\tag $\Cal W$
$$
 
The complexes on the axes are the usual Hochschild complexes. The
$n$'th column is the complex $\ap\pot{\g B}C_{\star}(\g B,Q\pot{}\Cal
A^{\pot{}(n-1)})$ where $C_{\star}(\g B,Q\pot{}\g A^{\pot{}(n-1)})$ is
the normalized bar resolution of the left $\g B$-module $Q\pot{}\Cal
A^{\pot{}(n-1)}$. Similarly, the $m$'th row is the complex $\bq\pot{\Cal
A}C_{\star}(\g A,P\pot{}\g B^{\pot{}(m-1)})$. For details, see
[G2]. Concerning exactness we have 

\proclaim{3.4 Lemma} Let ($\Cal W^*$) be the dual double co-complex of
($\Cal W$) and suppose that $\g A\simm\g B$. Then there is vertical
exactness at places $(n,i)$ for $i=0,1$ and $n\geq1$  and horizontal
exactness at places $(i,n)$ for $i=0,1$ and $n\geq1$ in $(\Cal W^*$. If $\g A$
has a BLAI, then columns of $(\Cal W^*)$ is acyclic except possibly on
the vertical edge. 
\endproclaim  
\demo{Proof} With minor modifications the proofs of [G2,Lemma 3.1] and
[G3,Theorem 4.6] can be adapted to the present situation, so the reader
is referred to these references.
\enddemo 

Applying this to algebras of approximable operators we get 

\proclaim{3.5 Theorem} Suppose that $\A X$ and $\A Y$ are self-induced
and that $\A X\simm\A Y$. If $X$ has the bounded approximation
property, then there is an injection
$$
\Cal H^1(\A Y,\A Y^*)\to \Cal H^1(\A X,\A X^*).
$$
In particular, if $X$ has the bounded approximation property and $\A X$ is
\wae, then $\A Y$ is \wae.
\endproclaim
\demo{Proof} According to Lemma 3.4, the double co-complex $(\Cal
W^*)$ corresponding to $\g A=\A X,\;\g B=\A Y$ satisfies the
hypothesis of Lemma 3.1 to conclude injectivity.
\enddemo
\remark{3.6 Remark} From [G3, Corollary 4.9] it follows that if $X$ and $Y$ both have
the bounded approximation property, then $\A X\simm\A Y$ implies 
$$
\Cal H^n(\A Y,\A Y^*)\cong \Cal H^n(\A X,\A X^*)\text{ for all }n\in\Bbb N.
$$
\endremark

\subheading{4. Some illustrative applications}

As mentioned in the introduction our approach will be to establish
\way\ for some key examples and then conclude \way\ for other Banach
algebras by means of the relation $\simm$. Towards this end we start
by
\proclaim{4.1 Theorem} Let $X$ be an arbitrary Banach space. Then for any
$1\leq p<\infty$
$$
\Cal H^1(\Bn{\ell_p(X)},\A{\ell_p(X)^*}=\set0\;.
$$ 
\endproclaim 
\demo{Proof} Since $\cl{\A{\ell_p(X)}^2}=\A{\ell_p(X)}$ a derivation
$D\colon\Bn{\ell_p(x)}\to\A{\ell_p(X)}^*$ is given by its restriction
to $\A{\ell_p(X)}$, that is, there is $T\in\Bn{\ell_p(X)^*}$ such that
$$
\l A,D(S)\r=\tr((SA-AS)^*T),\quad A\in\F{\lp X},S\in\Bn{\lp X}\;.
$$
We want to find $\lambda\in\Bbb C$ such that $T-\lambda\bold1$ is
integral. We start by seting up some notation. We view $\Bn{\lp X}$ as
consisting of infinite matrices with each entry an operator from $\Bn
X$. For $n\in\Bbb N$ we let $V_n$ and $H_n$ denote the left and right
shifts by $n$ places. We let $\g M=\{\bold W\in\F{\lp X}\mid\exists n\in\Bbb
N\colon \bold WH_n=V_n\bold W=0\}$, i.e\. $\g M$ is the dense
subalgebra of $\A{\lp X}$ consisting of matrices of the form
$$
\bold W=\2x2 W000\;,
$$
where $W$ denotes a finite square matrix with entries from $\F
X$ and the 0's represent infinite 0-matrices of the appropriate
size. Let $\bold W\in\g M$ and choose a $(d\times d)$-matrix $W$ to
represent $\bold W$. For $N\in\Bbb N$ we write $\Delta_{N,d}(W)$ for the
matrix obtained by repeating the matrix  $W$ 
along the diagonal $N$ times, i.e\.
$$
\Delta_{N,d}(W)=\pmatrix W&&\\
{}&\ddots&\\
{}&&W\\
\endpmatrix.
$$ 
Note that 
$$
\nrm{\2x2{\Delta_{N,d}(W)}000}=\nrm{\bold W}.
$$ 
Note also that a given $\bold
W\in\g M$ can be represented by different matrices $W$, since we may
add 0-rows and 0-columns. In order to prove that $T-\lambda\bold1$ is
integral it suffices to prove that there is a constant $C>0$ such that
$$
\vert\tr(\bold W^*(T-\lambda\bold1)) \vert\leq C\nrm{\bold W},\quad
\bold W\in\g M.
$$

Let $\bold W\in\g M$. Then 
$$
\vert\tr((\bold W-H_n\bold WV_n)^*T)\vert=\vert\l H_n\bold W,D(V_n)\r
\vert\leq\nrm{\bold W}\nrm D\;.
$$
It follows that the sequence $(\tr((H_n\bold WV_n)^*T))$ is bounded. Let
$\LIM$ be a Banach limit and define a linear functional (possibly unbounded)
$f\colon\g M\to\Bbb C$ by
$$
f(\bold W)=\LIM(\tr((H_n\bold WV_n)^*T)),\quad \bold W\in\g M\;.
$$
We now prove that $f(\bold U\bold W)=f(\bold W\bold U),\;\bold U,\bold W\in\g M$. By including some
0-entries, if necessary, we may suppose that $\bold U$ and $\bold W$ are
represented by matrices, $U$ and $V$ respectively, of equal size, say $d\times d$. Let for
$n\in\Bbb N$ the  0-matrix of size $n\times n$ be denoted $0_n$ and consider
the operators in $\g M$ given by the matrices
$$
\align
R_n(W)&=\pmatrix{0_n}&0&0&0\\
        0&0_{Nd}&\Delta_{N,d}(W)&0 \\
        0&0&0&0\endpmatrix\\
S_n(U)&=\pmatrix{0_n}&0&0\\
        0&0_{Nd}&0\\
        0&\Delta_{N,d}(U)&0\\
        0&0&0\endpmatrix\;. 
\endalign
$$
Then $\nrm{R_n(W)}=\nrm{\bold W}$ and $\nrm{S_n(U)}=\nrm{\bold U}$. From the identity
$$
\l R_n(W),D(S_n(U)\r=\tr\left(\pmatrix 0_n&0&0&0\\
                                       0&\Delta_{N,d}(WU)&0&0\\
                                       0&0&-\Delta_{N,d}(UW)&0\\
                                       0&0&0&0\endpmatrix^*T\right)
$$
and from translation invariance of the Banach limit we get
$$
\align
\vert f(\bold W\bold U-\bold U\bold W)\vert&=\frac1N\vert\LIM\l R_n(W),D(S_n(U)\r\vert\\
                  &\leq\frac1N\nrm D\nrm{\bold W}\nrm{\bold U}\;.      
\endalign
$$
Since $N$ is arbtrary, we arrive at $f(\bold W\bold U)=f(\bold U\bold W)$. It follows that
there is $\lambda\in\Bbb C$ such that $f(\bold W)=\lambda\tr(\bold W)$. This is the
$\lambda$ we are looking for:
$$
\align
&\tr\left(\2x2W000^*(T-\lambda\bold1)\right)=\\
&                       \tr\left(\pmatrix W&0&0&0\\
                                0&{0_n}&0&0\\
                                0&0&-W&0\\
                                0&0&0&0\endpmatrix^*(T-\lambda\bold1)\right)+
                        \tr\left(\pmatrix0&0&0&0\\
                                0&{0_n}&0&0\\
                                0&0&W&0\\
                                0&0&0&0\endpmatrix^*(T-\lambda\bold1)\right)=\\
&                       \tr\left(\pmatrix W&0&0&0\\
                                0&{0_n}&0&0\\
                                0&0&-W&0\\
                                0&0&0&0\endpmatrix^*T\right)+
                        \tr\left(\pmatrix0&0&0&0\\
                                0&{0_n}&0&0\\
                                0&0&W&0\\
                                0&0&0&0\endpmatrix^*T\right)-\lambda\tr(W).
\endalign
$$
Since
$$
\pmatrix W&0&0&0\\
 0&{0_n}&0&0\\
0&0&-W&0\\
0&0&0&0\endpmatrix =\left[\2x2 0W00;\2x200P0\right]
$$
for an appropriate coordinate projection $P$, we get by taking $\LIM$
that
$$
\vert\tr(\bold W^*(T-\lambda\bold1)\vert=\vert\tr\left(\2x2
W000^*(T-\lambda\bold1)\right)\vert\leq\nrm D\nrm{\bold W}\;, 
$$
which is want we wanted.
\enddemo
\proclaim{4.2 Corollary} Let $X$ be a Banach space. Then
$$
\A{\ell_p(X)} \text{ is \wae}\iff\A{\ell_p(X)}\text{ is self-induced.}
$$
\endproclaim
\demo{Proof} $\ell_p(X)$ is decomposable, so if $\A{\ell_p(X)}$ is
\wae, then it is self-induced by Theorem 3.8. If $\A{\ell_p(X)}$ is
self-induced, then every derivation
$$
D\colon\A{\ell_p(X)}\to\A{\ell_p(X)}^*
$$ 
can be extended to a derivation $\tilde
D\colon\Bn{\ell_p(X)}\to\A{\ell_p(X)}^*$, which by Theorem 4.1 is inner.
\enddemo
\remark{4.3 Remark} In the proof of Theorem 4.1 the only properties of
$\ell_p(X)$ we used were (i): there is a constant $C>0$ with
$\nrm{\Delta_{N,d}(W)}\leq C\nrm W$ for all $W$, (ii): $\g M$ is dense
in $\A{\ell_p(X)}$. The latter is equivalent to $\lim_n H_n A V_n=0$
for all $A\in\A{\ell_p(X)}$. Hence there are many other Banach spaces of
sequences from $X$ for which the proof works, notably
$c_0(X)$. However, in the present paper we shall only make use of the
spaces $\ell_p(X),\;1\leq p<\infty$.
\endremark 

\smallskip
The next result concerning \way\ of $\A{\LP X}$ strengthens Theorem 4.1
of [B2] by weakening the hypothesis '$X^*$ has the bounded
approximation property' to '$X$ has the bounded approximation
property'. The data of the space $\LP X$ are a measure space
$(\Omega,\Sigma,\mu)$ and a sequence $(\Omega_n)$ of pairwise
disjoint sets in $\Sigma$ with $0<\mu(\Omega_n)<\infty$ (to avoid
simply dealing with the case $X\oplus\cdots\oplus X$). Without loss of
generality we may further assume that $\mu$ is a probability measure,
since every compact set in $\LP X$ has $\sigma$-finite support.

\proclaim{4.4 Theorem} Let $X$ be a Banach space with the bounded
approximation property. Then $\A{\LP X}$ is \wae. 
\endproclaim

\demo{Proof} Since $X$ has the bounded approximation property, the
same is true for the spaces $\lp X$ and $\LP X$. In particular $\A{\lp
X}$ is self-induced and therefore \wae\ by Corollary 4.2. Thus we may
prove the theorem by showing $\A{\LP X}\simm\A{\lp X}$. First we give
some notation and well-known facts.  A {\it mesh} $\g m=\set{E_n\mid
n\in\Bbb N}$ is a partition 
$\Omega=\bigcup_{n=1}^\infty E_n$ into pairvise disjoint
measurable sets. A mesh $\g m$
defines a norm-1 projection $P_{\g m}\in\Bn{\LP}$ by the rule
$$
P_{\g m}(f)=\sum_{E\in\g m,\mu(E)\neq0}(\frac1{\mu(E)}\int_{E}f\,d\mu)\xi_{E}
$$
The set $\set{\text{meshes}}$ is  ordered by refinement and $\lim_{\g
m\to\infty}P_{\g m}=\bold1_{\A{\LP X}}$ uniformly on compacta.
For a mesh-projection
$P_{\g m}$ the range is isometrically isomorphic to $\ell_p^\kappa(X)$,
where $\kappa$ is the cardinality of $\set{E\in\g m\mid \mu(E)>0}$. In
particular $\LP X$ has a complemented subspace isometric to $\lp X$,
so that 
$$
\LP X\cong\LP X\oplus\lp X.
$$
Since $\A{\lp X}$ has a bounded left approximate identity, it follows
that $\A{\lp X}$ factors approximately through $\A{\LP X}$. To show
that $\A{\LP X}$ factors approximately through $\A{\lp X}$ first note
that $\set{P_{\g m}G\mid G\in\A{\LP X}_1, \g m\text{ a mesh}}$ is
dense in $\A{\LP X}_1$. It follows that each $A\in\A{\LP X}_1$ is the
sum of a series
$$
A=\sum_1^\infty 2^{-n}P_nG_n,
$$
where the $P_n$'s are mesh-projections and $\set{G_n}\subseteq \A{\LP
X}_1$. Identifying the ranges of mesh-projections with the appropriate
$\lp X$-spaces we get  
$$
A=\sum_1^\infty 2^{-n}\tilde P_n\tilde G_n,
$$
with $\tilde P_n\in\Bb{\lp X}{\LP X}_1,\tilde G_n\in\Aa{\LP X}{\lp
X}_1$. Since $\A{\LP X}$ has a bounded left approximate identity, we
conclude that $\A{\LP X}$ factors approximately through $\A{\lp
X}$. Alltogether $\A{\LP X}\simm\A{\lp X}$.
\enddemo 

\smallskip
In order to facilitate the use of factorization properties, the
generalization given by Blanco in [B2] of Johnson's $C_p$-spaces is
very useful. We quote it here:
\definition{4.5 Definition} A Banach space $J$ is called a {\it
Johnson space} if it has the form $(\oplus_1^\infty G_n)_p, p=0,\;1\leq
p<\infty$, where $(G_n)_n$ is a sequence of finite-dimensional Banach
spaces such that for each $i\in\Bbb N$ the set $\{n\in\Bbb N\mid
G_n\cong X_i\text{ isometrically}\}$ is infinite.

Let $J=(\oplus G_n)_p$ be a Johnson space. A Banach space $X$ is called
a {\it J-space} if there is $\lambda\geq1$ such that for every
finite-dimensional  subspace $E$ of $X$, there is a subspace $G$ of $X$
containing $E$ such that the Banach-Mazur distance
$d(G,G_i)\leq\lambda$ for some $i\in\Bbb N$. 
\enddefinition

The usefulness of these notions lies in

\proclaim{4.6 Proposition} Let $J=(\oplus_1^\infty
G_n)_p$ be a Johnson space, and let $X$ be a $J$-space. Then
$\A J$ is \wae, and $\A X$ factors approximately through $\A J$.
\endproclaim
\demo{Proof} It is an immediate consequence of Corollary 4.2 that $\A
J$ is \wae. Let $A\in\F X$ and choose
$\operatorname{range}(A)\subseteq G$ and corresponding $G_i$ in accordance with the definition of $X$
being a $J$-space. This gives a factorization
$$
\CD G_i&@>\iota>>& J\\
@AVAA&&@VV U V\\
X&@>>A>&X
\endCD
$$ 
with $\Vert U\Vert\Vert V\Vert\leq\lambda\nrm A$. The claim now follows from
Lemma 2.6.
\enddemo

The next result is Proposition 3.3 of [B2]. 

\proclaim{4.7 Proposition {\rm [Blanco]}} Let $J=(\oplus_1^\infty
G_n)_p$ be a Johnson space, and let $X$ be a $J$-space. Then
$\A{X\oplus J}$ is \wae.
\endproclaim
\demo{Proof} Since $J$ has the bounded approximation property we
obviously have that $\A J$ factors through $\A{X\oplus J}$ and from
Proposition 4.2 it follows that $X\oplus J$ factors approximately through
$J$, i.e. $\A{X\oplus J}\simm\A J$. Since $\A J$ is \wae\ by Corollary
4.2, the proof is concluded by an appeal to Theorem 3.5.
\enddemo

As a final illustration, we shall look at the James spaces $\g
J_p$. Blanco [B2] shows that $\A{\g J_p}$ is \wae, by showing that
there is a Johnson space $J_p$ such that $\g J_p$ is a $J_p$-space and
$\g J_p\cong\g J_p\oplus J_p$, whence the result follows from
Proposition 4.7. Using the relation $\simm$ makes it possible to
extend this result to vector-valued James spaces. We start by briefly
recalling basic properties of the spaces $\g J_p$ following the
notation of [B2]. Let $1<p<\infty$ and let $(\alpha_n)\in\Bbb C^{\Bbb
N}$. Define $\jnrm\cdot$ by

$$
\jnrm{(\alpha_n)}=\sup\set{(\sum_{n=1}^{m-1}\vert\alpha_{i_n}-\alpha_{i_{n+1}}\vert^p)^{\frac1p}\mid
i_1<\dots<i_m,\;m\geq2}.
$$
Then
$$
\g J_p=\set{(\alpha_n)\in\Bbb C^{\Bbb
N}\mid\jnrm{(\alpha_n)}<\infty,\lim_n\alpha_n=0}.
$$

With this norm $\g J_p$ is a Banach space. The sequences
$e_n=(\delta_{kn})_k$ form a normalized, 1-unconditional basis, $\bold
e$, for $\g J_p$, the {\it canonical basis\/}. We now define vector
valued James spaces. But we start with a general setting which is a
special case of the spaces described in [Lau].

\definition{4.8 Definition} Let $E$ be a Banach space with a
normalized 1-unconditional basis $\bold b=\set{b_1,\dots}$, and let $X$
be any Banach space. 
Then we define $E\wot bX$ as 
$$ 
E\wot bX=\set{(x_n)\in X^{\Bbb N}\mid \sum_1^\infty\nrm{x_n}b_n\in E}
$$
with the norm $\nrm{(x_n)}_{\bold b}=\nrm{\sum_1^\infty\nrm{x_n}b_n}$.
The {\it $X$-valued James $p$-space} is $\g J_p(X)=\g J_p\wot eX$,
where $\bold e$ is the canonical basis. In
this case we use the notation $\jnrmx\cdot=\nrm\cdot_{\bold e}$. 
\enddefinition

It is straightforward to verify that $(E\wot bX,\nrm\cdot_{\bold b})$ is a
Banach space, which may be viewed as a completion of the algebraic
tensor product $E\otimes X$, if we identify
$(\sum_1^\infty\alpha_nb_n)\otimes x\in E\otimes X$ with
$(\alpha_nx)\in E\wot bX$. Note that in this picture $\nrm\cdot_{\bold
b}$ is a cross-norm: $\nrm{(\sum_1^\infty\alpha_nb_n)\otimes x}_{\bold
b}=\nrm{\sum_1^\infty\alpha_nb_n}\nrm x$. In accordance with this
picture consider $S\in \Bn E,\;T\in \Bn X$. If the linear map
$S\otimes T\colon E\otimes X\to E\otimes X$ extends to a bounded
operator $E\wot bX\to E\wot bX$, the latter will be denoted $S\wot
bT$. 

\medskip
The result above by Blanco is the case $X=\Bbb C$ of the following
theorem.

\proclaim{4.9 Theorem} Let $X$ be a Banach space. If $X$ has the
bounded approximation property, then the Banach algebra $\A{\g J_p(X)}$
is \wae\ for every $1<p<\infty$. 
\endproclaim
\demo{Proof} For each $n\in\Bbb N$, define a closed subspace of $\g
J_p(X)$ by $\g
J_{p,n}(X)=\set{x_1e_1+\dots+x_ne_n\mid x_1,\dots,x_n\in X}$. Let $(G_k)$ be a sequence of
Banach spaces obtained by repeating each $\g
J_{p,n}(X)$ infinitely many times. Define  a Banach space by $J_p(X)=(\oplus
G_k)_p$. We prove that $\A{\g J_p(X)}\simm\A{J_p(X)}$. Invoking Theorem
3.5 and  Corollary 4.2, the claim
follows, since both spaces $\g J_p(X),\;J_p(X)$ have the bounded
approximation property and $\lp{J_p(X)}\cong J_p(X)$. 
Let $P_n\colon\g J_p(X)\to\g J_p(X)$ be the canonical projection
onto $\g J_{p,n}(X)$. Then $P_nA@>>>A$ for all $A\in\A{\g J_p(X)}$. Each
$P_n$ having an obvious factorization
$$
P_n=\iota_nQ_n,\quad Q_n\in\Aa{\g J_p(X)}{J_p(X)},\iota_n\in\Aa{J_p(X)}{\g
J_p(X)},\;\nrm{Q_n}=\nrm{\iota_n}=1,
$$
it follows that each $A\in\A{\g J_p(X)}$ has a decomposition
$$
A=\sum_1^\infty T_nS_n,\quad S_n\in\Bb{\g J_p(X)}{J_p(X)},T_n\in\Bb{J_p(X)}{\g
J_p(X)},
$$
with $\sum\nrm{S_n}\nrm{T_n}\leq2\nrm A$. Since $\g J_p(X)$ has the
bounded approximation property, we may write $A=A_1A_2A_3,\;
A_1,A_2,A_3\in\A{\g J_p(X)}$. It follows that $\A{\g J_p(X)}$ factors
approximately through $\A{J_p(X)}$.

To prove that $\A{J_p(X)}$ factors approximately through $\A{\g
J_p(X)}$ we just note that Blanco's decomposition of $\g J_p$ works
equally well for $\g J_p(X)$ with the same proof, so that we have $\g
J_p(X)\cong J_p(X)\oplus J_p(X)$. Since $J_p(X)$ has the bounded
approximation property, we may factor as desired.
\enddemo

\subheading{5. Conclusion} As demonstrated, many questions of \way\ of
Banach algebras (notably of the type $\A X$) can be approached using
factorization of Morita equivalence type. This has been illustrated by
giving a framework behind much of the reasoning in Blanco's papers
[B1] and [B2]. We would like to raise some questions related to this. 

\remark{5.1 Question} In Proposition 4.2 of [B1] Blanco shows that if
$P$ is a Banach space such that $P$ and $P^*$ both have cotype 2, then
$\A{\ell_2(P)}$ is \wae. Part of his argument consists in using a
factorization theorem by Pisier (Theorem 4.1 of [P1]), which combined
with Lemma 2.6 shows that $\A{\ell_2(P)}$ is self-induced. However, $P$
can be chosen such that $\A P$ is not \wae.

In the same paper Blanco constructs a reflexive space $E$ with an
unconditional basis such that $\A E$ is not \wae. Since $E$ has the
bounded approximation property, $\A{\ell_p(E)}$ is self-induced and
hence \wae\ for $1\leq p<\infty$. 

The spaces $\ell_p(X)$ have the form $\ell_p\wot eX$ and are {\it
tight tensor products} in the sense of [GJW]. Thus we may view
$\A{\ell_p(X)}$ as a tensor product $\A{\ell_p}\wt\otimes\A X$. 
The preceeding paragraphs can be phrased as a stabilizing effect of
the functor $\A{\ell_p}\wt\otimes-$, in liking with stabilizing in
the theory of $C^*$-algebras. This leads to
\medskip
 Is $\A{\ell_2(X)}$ self-induced for all $X$? Is $\A{\lp
X},\;1\leq p<\infty$? A test case would be the space constructed in [P2] for which
multiplication is not surjective.
\endremark
\remark{5.2 Question} Our reasoning has relied on matrix-structures
with a certain uniformity loosely speaking enabling us to shift
matrices around. Hence the Tsirelson space $T$, for which Blanco
established \way\ of $\A T$, presents a possible shortcoming. But in
proving that $\A T$ is \way, it suffices to prove that $\A X$ is \wae\
for some 'nice' $X$ such that $\A T\simm\A X$. This leads to 
\medskip
What are the spaces $X$ such that $\A X$ and $\A T$ are Morita
equivalent? Which among these have an unconditional basis? 

\endremark

\enddocument